\documentclass[11pt]{article}
\usepackage{amssymb}
\usepackage{graphicx}
\usepackage{setspace}
  \usepackage{paralist}
  \usepackage{longtable}
   \usepackage{multirow}
    \usepackage{rotating}
\usepackage{marginnote}
\usepackage{mathrsfs}
\usepackage{amsmath, amsthm, amssymb}
\usepackage{authblk}
\usepackage{graphicx}
\usepackage{float}
\usepackage{hyperref}
\usepackage[margin=1in]{geometry}
\usepackage{comment}
\usepackage{color}
\usepackage{cite}
\usepackage{indentfirst}

\allowdisplaybreaks[4]
\numberwithin{equation}{section}
\date{}
\newtheorem{theorem}{Theorem}[section]
\newtheorem*{theorem*}{Theorem}

\newtheorem{corollary}[theorem]{Corollary}

\theoremstyle{remark}

\begin{document}
\renewcommand{\thefootnote}{\fnsymbol{footnote}}

\footnotetext{1. School of Mathematical Sciences, CMA-Shanghai, Shanghai Jiao Tong University, China( xmq157@sjtu.edu.cn);}
\footnotetext{2. Partially supported by NSFC-12031012 and NSFC-11831003.}
\renewcommand{\thefootnote}{\arabic{footnote}}

\title{On super polyharmonic property of high-order fractional Laplacian}
\author{Meiqing Xu}

\maketitle
\begin{abstract}
\noindent Let $0<\alpha<2$, $p\geq 1$, $m\in\mathbb{N}_+$. Consider $u$ to be the positive solution of the PDE
\begin{equation}\label{abstract PDE}
     (-\Delta)^{\frac{\alpha}{2}+m} u(x)=u^p(x) \quad\text{in }\mathbb{R}^n.
\end{equation}  In \cite{cao2021super}(Transactions of the American mathematical society, 2021), Cao, Dai and Qin showed that, under the condition $u\in\mathcal{L}_\alpha$, (\ref{abstract PDE}) possesses super polyharmonic property $(-\Delta)^{k+\frac{\alpha}{2}}u\geq 0$ for $k=0,1,...,m-1$. In this paper, we show another kind of super polyharmonic property $(-\Delta)^k u> 0$ for $k=1,...m$ under different conditions $(-\Delta)^mu\in\mathcal{L}_\alpha$ and $(-\Delta)^m u\geq 0$. Both kinds of super polyharmonic properties can lead to the equivalence between (\ref{abstract PDE}) and the integral equation $u(x)=\int_{\mathbb{R}^n}\frac{u^p(y)}{|x-y|^{n-2m-\alpha}}dy$. One can classify solutions to (\ref{abstract PDE}) following the work of \cite{chen2003qualitative}\cite{chen2} by Chen, Li, Ou.


\end{abstract}

\section{Introduction}
In this paper we consider the partial differential equation
\begin{align}\label{main equation}
  (-\Delta)^{\frac{\alpha}{2}+m}u(x)=u^p(x)\quad\text{in } \mathbb{R}^n,
\end{align}
under conditions that $n\geq 2$, $0<\alpha<2$, $p\geq 1$, $m\in\mathbb{N}_+$.

For $\alpha$ taking any real number between 0 and 2, $(-\Delta)^\frac{\alpha}{2}$ is a nonlocal differential operator defined by
\begin{equation}\label{fraction definition}
  (-\Delta)^\frac{\alpha}{2}u(x)=C_{n,\alpha}PV \int_{\mathbb{R}^n} \frac{u(x)-u(y)}{|x-y|^{n+\alpha}}dy.
\end{equation}
where PV is the Cauchy principal value. (\ref{fraction definition}) is valid for $u\in C_{loc}^{[\alpha],\{\alpha\}+\epsilon}(\mathbb{R}^n) \cap\mathcal{L}_\alpha(\mathbb{R}^n)$, where
\begin{equation}
  \mathcal{L}_\alpha(\mathbb{R}^n)=\{u:\mathbb{R}^n\rightarrow\mathbb{R}| \int_{\mathbb{R}^n} \frac{|u(x)|}{(1+|x|^{n+\alpha})}dx <\infty\}.
\end{equation}
The super polyharmonic property $(-\Delta)^{k+\frac{\alpha}{2}}u(x)\geq 0$ for $k=0,1,...,m-1$ of (\ref{main equation}) was first proved in \cite{cao2021super}(Transactions of the American mathematical society, 2021) by Cao, Dai and Qin under the condition that $u\in\mathcal{L}_\alpha$. In fact they worked on more general equations of the form $(-\Delta)^{m+\frac{\alpha}{2}}u(x)=f(x,u,Du,...)$. Their idea is to use Green's formula and Poisson kernel to represent the PDE.

Inspired by their method, here in this paper we obtain another kind of super polyharmonic property under a quite different definiton. We write $(-\Delta)^{\frac{\alpha}{2}+m}u(x)$ as $(-\Delta)^\frac{\alpha}{2}(-\Delta)^mu(x)$, and therefore it is nature to require $(-\Delta)^mu(x)\in\mathcal{L}_\alpha$. Moreover we require $(-\Delta)^mu(x)\geq 0$. And the result is a different kind of super polyharmonic property $(-\Delta)^ku(x)\geq0$ for $k=1,...,m$.
\begin{theorem}\label{thm super harmonic}
  Let $n\geq 2$, $0<\alpha<2$, $p\geq 1$, $m\in\mathbb{N}_+$. Suppose $(-\Delta)^m u\in\mathcal{L}_\alpha(\mathbb{R}^n)$, $u\in C^\infty(\mathbb{R}^n)$ and $u$ is a positive solution of
(\ref{main equation}).
  Moreover, suppose $(-\Delta)^m u\geq 0$ in $\mathbb{R}^n$. Then for $k=1,...,m$,
  \begin{equation}
    (-\Delta)^k u(x)> 0\quad\text{in } \mathbb{R}^n.
  \end{equation}
\end{theorem}
The proof for $(-\Delta)^mu\geq0$ stays unsolved and will be worked on in the future.

In theorem (\ref{thm super harmonic}) we require $u\in C^\infty(\mathbb{R}^n)$, but it is not necessary. In fact one only needs $C_{loc}^{2m+[\alpha],\{\alpha\}+\epsilon}(\mathbb{R}^n)$(see \cite{silvestre}, proposition 2.4) to make sure $(-\Delta)^{\frac{\alpha}{2}+m} u$ is a continuous functions and its value is given by (\ref{fraction definition}). Nevertheless, $u\in C^\infty(\mathbb{R}^n)$ is enough since the precise regularity is not the main goal in this paper.

(\ref{main equation}) is closely related to the integral equation
\begin{equation}\label{integral}
  u(x)=\int_{\mathbb{R}^n}\frac{u^p(y)}{|x-y|^{n-2m-\alpha}}dy.
\end{equation}
Once super polyharmonic property holds, then it is quick to obtain the equivalence between (\ref{integral}) and (\ref{main equation}). Then one can use Chen, Li, Ou's work \cite{chen2003qualitative}\cite{chen4} about classification of solutions to integral equations and therefore obtain the classification of high-order Laplacian.
\begin{theorem}\label{thm equivalence}
  Suppose conditions in theorem \ref{thm super harmonic} hold, and
   \begin{equation}
       \int_{\mathbb{R}^n}\frac{u^p(y)}{|x-y|^{n-2m-\alpha}}dy<+\infty\quad\mbox{for any }x\in\mathbb{R}^n.
   \end{equation}
   Let $u$ be a positive solution of (\ref{main equation}).  Then $u$ also solves (\ref{integral}), and vice versa.
\end{theorem}

\begin{corollary}\label{classification}
   The solution $u$ as given in theorem \ref{thm equivalence} satisfies
   \begin{enumerate}
       \item For $0<2m+\alpha<n$, in the critical case  $p=\frac{n+2m+\alpha}{n-2m-\alpha}$, $u$ must have the form
       \begin{equation}
           u(x)=c(\frac{t}{t^2+|x-x_0|^2})^\frac{n-2m-\alpha}{2}
       \end{equation}
       with some constant $c=c(n,m,\alpha)$ and for some $t>0$ and $x_0\in\mathbb{R}^n$.
       \item  For $0<2m+\alpha<n$, in the subcritical case $1\leq p<\frac{n+2m+\alpha}{n-2m-\alpha}$, $u$ does not exist.
   \end{enumerate}
\end{corollary}
Corollary \ref{classification} is the immediate consequence of theorem \ref{thm equivalence} and classification results in \cite{chen2003qualitative}\cite{chen4}.

Super polyharmonic property of high-order Laplacian has long been studied. It leads directly to Liouville theorem and equivalence between integral equations like (\ref{integral}) and PDEs like (\ref{main equation}). For integer high-order Laplacian, Wei and Xu \cite{wei1999classification} first gave the super polyharmonic property for even order equations. After a dozen years Chen and Li \cite{chen2013super} proved the general case and covered Wei and Xu's work. There are also much work on other types of super polyharmonic property, see \cite{fang2012liouville}\cite{chen2013super2}\cite{dai2021classification}\cite{cheng2016liouville}\cite{zhuo2019liouville}.

For more results about fractional equation $(-\Delta)^\frac{\alpha}{2}=u^p$, polyharmonic equations and maximum principle about fractional Laplacian, please refer to \cite{cheng2017direct}\cite{zhuo2022classification}\cite{zhuo2017weighted}\cite{li2020non}.

\section{Proof of theorem \ref{thm super harmonic}}
\begin{proof}
  Denote $(-\Delta)^k u(x)$ as $u_k(x)$ for $k=1,...,m$. Then (\ref{main equation}) can be rewritten as
  \begin{equation}
  \begin{cases}
        -\Delta u=u_1\\
    -\Delta u_1=u_2\\
    ...\\
  -\Delta u_{m-1}=u_m\\
    (-\Delta)^\frac{\alpha}{2}u_m=u^p
    \end{cases}\quad\text{in } \mathbb{R}^n.
  \end{equation}
First show that if $(-\Delta)^m u$ is nonnegative, then it must be positive. Indeed, if there exists $x_0$ such that $(-\Delta)^m u(x_0)=u_m(x_0)=0$, then $u_m(x_0)-u_m(y)\leq 0$ for any $y\in\mathbb{R}^n$. So
\begin{equation}
  (-\Delta)^\frac{\alpha}{2}u_m(x_0)=\int_{\mathbb{R}^n}\frac{u_m(x_0)-u_m(y)}{|x-y|^{n+\alpha}} dy\leq 0,
\end{equation}
contradicted with $(-\Delta)^\frac{\alpha}{2}u_m(x_0)=u^p(x_0)>0$.

Now prove $u_k(x)>0$ for $k=1,...,m$ by contradiction. Suppose NOT, then there exists a largest integer $k\in\{1,2,...,m-1\}$ such that $u_k$ is less or equal to $0$ somewhere. Without loss of generality suppose $u_k(0)\leq 0$. Then there are two cases:
\begin{enumerate}
  \item\label{=0 case} $u_k(0)=0$, and $u_k\geq 0$ in $\mathbb{R}^n$.
  \item\label{negative case} $u_k(0)<0$.
\end{enumerate}
Case \ref{=0 case} is impossible. The reason is, under this case, $0$ is a minimum point of $u_k$. However $-\Delta u_k=u_{k+1}\geq 0$, by maximum principle $u_k$ cannot obtain minimum. So only need to consider case \ref{negative case}.

We argue by two steps to see such $k$ does not exist.

Step 1, $k$ cannot be even. Denote by $\overline{v}$ the spherical average of $v$ centered at 0, i.e.
\begin{equation}
  \overline{v}(r)=\frac{1}{|\partial B(0,r)|}\int_{\partial B(0,r)}v(y) dS_y.
\end{equation} Use the well-known property $\Delta\overline{v}=\overline{\Delta v}$ and $-\Delta \overline{v}(r)=-\frac{1}{r^{n-1}}(r^{n-1}\overline{v}')'$ (here $'$ means taking derivative about $r$) to get
\begin{align}
  -\Delta \overline{u_k}(r) =-\frac{1}{r^{n-1}} (r^{n-1}\overline{u_k}')'= \overline{u_{k+1}}(r)>0\quad\mbox{for any }r>0.
\end{align}
Then $(r^{n-1}\overline{u_k}')'<0$. Integrate on both sides to get $r^{n-1}\overline{u_k}'<0$. So
\begin{align}
  \overline{u_k}'(r)<0\quad\mbox{for any }r>0.
\end{align}
Thus
\begin{equation}
  \overline{u_k}(r)<\overline{u_k}(0)=-a_k<0\quad\text{for any }r>0.
\end{equation}
Using $-\Delta \overline{u_{k-1}}=\overline{u_k}$ one obtains
\begin{equation}
    -\frac{1}{r^{n-1}}(r^{n-1} \overline{u_{k-1}}')'<-a_k\quad\mbox{for any }r>0,
\end{equation}
and in turn implies
\begin{equation}
    \overline{u_{k-1}}(r)>\overline{u_{k-1}}(0)+\frac{a_k}{2n}r^2\quad\mbox{for any }r>0.
\end{equation}
If $k$ is even, repeat the process above to get
\begin{equation}
  \overline{u}(r)<\overline{u}(0)+c_1r^2-c_2r^4+...-c_{k}r^{2k}\quad\text{for any }r>0,
\end{equation}where we know $c_{k}>0$ but do not know the sign of other $c_i$, $i=1,...,k-1$. For $r$ sufficiently large, $\overline{u}(r)$ attains negative value, contradicted to $u>0$.

Step 2, $k$ cannot be odd. If $k$ is odd, repeat the similar process above to get
\begin{equation}
  \overline{u}(r)>\overline{u}(0)-c_1r^2+c_2r^4-...+c_kr^{2k}\quad\text{for any }r>0,
\end{equation}
where we know $c_k>0$ but do not know the sign of other $c_i$, $i=1,...,k-1$. Hence there exists $r_0>0$ and a positive constant $M$ such that
\begin{equation}
    \overline{u}(r)>M>0,\quad\mbox{for any }r>r_0.
\end{equation}
Let $u_m=v_1+v_2$ with $v_1,v_2$ satisfying
\begin{equation}
  \begin{cases}
    (-\Delta)^\frac{\alpha}{2}v_1=u^p, & \mbox{in } B(0,R), \\
    v_1=0, & \mbox{on }\mathbb{R}^n\backslash B(0,R),
  \end{cases}
\end{equation} and
\begin{equation}
  \begin{cases}
 (-\Delta)^\frac{\alpha}{2}v_2=0, & \mbox{in } B(0,R), \\
    v_2=u_m, & \mbox{on }\mathbb{R}^n\backslash B(0,R),
  \end{cases}
\end{equation}
for any $R>0$.
By \cite{bucur} it follows
\begin{equation}
u_m(x)=\int_{B(0,R)}G^\alpha_R(x,y)u^p(y)dy+\int_{\mathbb{R}^n\backslash B(0,R)}P_R^\alpha(y,x)u_m(y)dy.
\end{equation}
Here $G^\alpha_R(x,y)$ and $P_R^\alpha(y,x)$ represent the Green's function and Poisson kernel of $\alpha$-Laplacian in $B(0,R)$, respectively. The formulas are

\begin{align}
  \label{fraction green} G^\alpha_R(x,y)=\frac{\kappa(n,\alpha)}{|y-x|^{n-\alpha}} &\int_{0}^{R_0(x,y)}\frac{t^{2\alpha-1}} {(t+1)^\frac{n}{2}}dt,\\
  P_R^\alpha(y,x)=c(n,\alpha) (\frac{R^2-|x|^2}{|y|^2-R^2}) ^\frac{\alpha}{2} \frac{1}{|x-y|^n}\quad&\mbox{for }x\in B(0,R) \mbox{ and }y\in\mathbb{R}^n\backslash \overline{B}(0,R)
\end{align}
with
\begin{equation}
  R_0(x,y)=\frac{(R^2-|x|^2)(R^2-|y|^2)}{R^2|x-y|^2}
\end{equation}
and $\kappa(n,\alpha),c(n,\alpha)$ are constants depending only on $\alpha$ and $n$. Then
\begin{equation}\label{geq0}
\begin{aligned}
u_m(0)&=\int_{B(0,R)}\frac{\kappa(n,\alpha)}{|y|^{n-\alpha}} (\int_{0}^{\frac{R^2}{|y|^2}-1} \frac{t^\frac{\alpha}{2}-1}{(1+t)^\frac{n}{2}}dt)u^p(y)dy + c(n,\alpha) \int_{\mathbb{R}^n\backslash B(0,R)} \frac{R^\alpha}{(|y|^2-R^2)^\frac{\alpha}{2}}\cdot \frac{u_m(y)}{|y|^n}dy\\
&=\kappa(n,\alpha)\int_0^R r^{\alpha-1}(\int_0^{\frac{R^2}{r^2}-1}\frac{t^{\frac{\alpha}{2}-1}}{(1+t)^\frac{n}{2}}dt)\overline{u^p} (r)dr+c(n,\alpha)\int_R^{+\infty}\frac{R^\alpha}{r(r^2-R^2)^\frac{\alpha}{2}}\overline{u_m}(r)dr\\
&\geq\kappa(n,\alpha)\int_0^R r^{\alpha-1}(\int_0^{\frac{R^2}{r^2}-1}\frac{t^{\frac{\alpha}{2}-1}}{(1+t)^\frac{n}{2}}dt)\overline{u}^p (r)dr.
\end{aligned}
\end{equation}

For $0<r<\frac{R}{2}$,
\begin{equation}
    \int_0^{\frac{R^2}{r^2}-1}\frac{t^{\frac{\alpha}{2}-1}}{(1+t)^\frac{n}{2}}dt\geq \int_0^3 \frac{t^{\frac{\alpha}{2}-1}}{(1+t)^\frac{n}{2}}dt,
\end{equation}
and the right hand side is a constant. Choose $R$ satisfying $R>2r_0$. Then
\begin{equation}\label{geq}
\begin{aligned}
\int_0^{r_0} r^{\alpha-1} (\int_0^{\frac{R^2}{r^2}-1}\frac{t^{\frac{\alpha}{2}-1}}{(1+t)^\frac{n}{2}}dt) \overline{u}^p (r)dr &\geq0,\\
\int_{r_0}^\frac{R}{2} r^{\alpha-1} (\int_0^{\frac{R^2}{r^2}-1}\frac{t^{\frac{\alpha}{2}-1}}{(1+t)^\frac{n}{2}}dt) \overline{u}^p (r)dr&\geq C\int_{r_0}^\frac{R}{2} r^{\alpha-1}dr\\
&= C_1R^\alpha-C_2\rightarrow +\infty\quad\mbox{as }R\rightarrow +\infty,\\
\int_\frac{R}{2}^R r^{\alpha-1} (\int_0^{\frac{R^2}{r^2}-1}\frac{t^{\frac{\alpha}{2}-1}}{(1+t)^\frac{n}{2}}dt) \overline{u}^p (r)dr&\geq 0.
\end{aligned}
\end{equation}
Combining (\ref{geq0}) and (\ref{geq}), it follows that $u_m(0)=\infty$, contradicted.
\end{proof}

\section{Proof of theorem \ref{thm equivalence}}
To prove the equivalence, one needs to show that any solution of one equation also solves another equation. It is easy that any solution of integral equation (\ref{integral}) also solves PDE (\ref{main equation}). So it suffices to show any solution of PDE solves integral equation.

Denote the Riesz potential as
\begin{equation}
    (I_\alpha f)(x)=\frac{1}{\gamma(\alpha,n)} \int_{\mathbb{R}^n} \frac{f(y)}{|x-y|^{n-\alpha}}dy,
\end{equation}
where
\begin{equation}
    \gamma(\alpha,n)=\frac{\pi^\frac{n}{2}2^\alpha\Gamma(\frac{\alpha}{2})}{\Gamma(\frac{n-\alpha}{2})}.
\end{equation}
The proof relies on two properties of Riesz potential( see \cite{stein2016singular}, chapter V, section 1.1).
\begin{enumerate}
     \item $(-\Delta)^\frac{\alpha}{2}(I_\alpha f)=f$ for $0<\alpha<n$.
    \item\label{prop riesz} $I_\alpha(I_\beta f)=I_{\alpha+\beta}f$ for $\alpha>0$, $\beta>0$ and $\alpha+\beta<n$.
\end{enumerate}
A brief outline is given here.
Rewrite (\ref{main equation}) as
  \begin{equation*}
  \begin{cases}
    (-\Delta)^\frac{\alpha}{2}u_m=u^p\\
      -\Delta u_{m-1}=u_m\\
       ...\\
           -\Delta u_1=u_2\\
              -\Delta u=u_1\\
    \end{cases}\quad\text{in } \mathbb{R}^n
  \end{equation*}
  and then prove by Liouville theorem and induction that
   \begin{equation*}
  \begin{cases}
  u_m(x)=f_1(x):= C\int_{\mathbb{R}^n}  \frac{u^p(y)}{|x-y|^{n-\alpha}}dy\\
     u_{m-1}(x)=f_2(x):=C\int_{\mathbb{R}^n}  \frac{f_1(y)}{|x-y|^{n-2}}dy\\
     ...\\
      u_1(x)=f_{m}(x):=C\int_{\mathbb{R}^n}  \frac{f_{m-1}(y)}{|x-y|^{n-2}}dy\\
    u(x)=C\int_{\mathbb{R}^n}  \frac{f_{m}(y)}{|x-y|^{n-2}}dy
    \end{cases}\quad\text{in } \mathbb{R}^n.
  \end{equation*}
  Finally by property \ref{prop riesz} of Riesz potential it follows that
  \begin{equation*}
      u(x)=C\int_{\mathbb{R}^n}  \frac{u^p(y)}{|x-y|^{n-2m-\alpha}}dy.
  \end{equation*}

\begin{proof}
Set $v_0^R(x)=\int_{B(0,R)}G_R^\alpha(x,y)u^p(y)dy$, where $G_R^\alpha(x,y)$ is the Green's function for fractional Laplacian as stated in (\ref{fraction green}). Then
\begin{equation}
    \begin{cases}
    (-\Delta)^\frac{\alpha}{2}v_0^R=u^p,&\quad\mbox{in } B(0,R),\\
    v_0^R=0,&\quad\mbox{on } \mathbb{R}^n\backslash B(0,R).\\
    \end{cases}
\end{equation}
Set $w_0^R=u_m-v_0^R$. It follows that
\begin{equation}
    \begin{cases}
    (-\Delta)^\frac{\alpha}{2}w_0^R=0,&\quad\mbox{in } B(0,R),\\
    w_0^R\geq 0,&\quad\mbox{on } \mathbb{R}^n\backslash B(0,R).\\
    \end{cases}
\end{equation}
By maximum principle( see \cite{silvestre}, proposition 2.17), $w_0^R\geq 0$. Sending $R\rightarrow\infty$, then $u_m(x)\geq C\int_{\mathbb{R}^n}\frac{u^p(y)}{|x-y|^{n-\alpha}}dy$.
Let $w_0=u_m-C\int_{\mathbb{R}^n}\frac{u^p(y)}{|x-y|^{n-\alpha}}dy$. Then
\begin{equation}
    \begin{cases}
    (-\Delta)^\frac{\alpha}{2}w_0=0,&\quad\mbox{in }\mathbb{R}^n,\\
    w_0\geq 0,&\quad\mbox{in }\mathbb{R}^n.\\
    \end{cases}
\end{equation}
By Liouville theorem( see \cite{zhuo2014liouville}, theorem 1), $w_0\equiv C_0\geq 0$, i.e
\begin{equation}\label{C0}
    u_m(x)=C\int_{\mathbb{R}^n}\frac{u^p(y)}{|x-y|^{n-\alpha}}dy+C_0.
\end{equation}
In fact $C_0$ can only be 0 thus $u_m(x)=C\int_{\mathbb{R}^n}\frac{u^p(y)}{|x-y|^{n-\alpha}}dy$, which will be proved later.

Let
\begin{align}
    f_1(x) &=C\int_{\mathbb{R}^n}\frac{u^p(y)}{|x-y|^{n-\alpha}}dy+C_0,\\
    v_1^R(x) &= \int_{\mathbb{R}^n}G^2_R(x,y)f_1(y)dy,\\
    w_1^R &=u_{m-1}-v_1^R.
\end{align}
where $G^2_R(x,y)$ is the Green's function of Laplacian $\Delta$ in $B(0,R)$. Then
\begin{equation}
    \begin{cases}
    -\Delta w_1^R=0 &\quad\mbox{in }B(0,R),\\
    w_1^R\geq 0, &\quad\mbox{on }\mathbb{R}^n\backslash B(0,R).
    \end{cases}
\end{equation}
By maximum principle of harmonic functions, $w_1^R\geq 0$. Sending $R\rightarrow\infty$, then $u_{m-1}(x)\geq C\int_{\mathbb{R}^n}\frac{f_1(y)}{|x-y|^{n-\alpha}}dy$. Let $w_1=u_{m-1}-C\int_{\mathbb{R}^n}\frac{f_1(y)}{|x-y|^{n-\alpha}|}dy$. Then
\begin{equation}
    \begin{cases}
    (-\Delta)^\frac{\alpha}{2}w_1=0,&\quad\mbox{in }\mathbb{R}^n,\\
    w_1\geq 0,&\quad\mbox{in }\mathbb{R}^n.\\
    \end{cases}
\end{equation}
By Liouville theorem by ones for harmonic functions, one obtains
\begin{equation}\label{C1}
    u_{m-1}(x)=C\int_{\mathbb{R}^n}\frac{f_1(y)}{|x-y|^{n-\alpha}}dy+C_1
\end{equation}
for some $C_1\geq 0$. Moreover, (\ref{C1}) implies $C_0=0$, otherwise
\begin{equation}
    +\infty>C\int_{\mathbb{R}^n}\frac{f_1(y)}{|x-y|^{n-\alpha}}dy\geq\int_{\mathbb{R}^n}\frac{C_0}{|x-y|^{n-\alpha}}dy=+\infty.
\end{equation}
Let
\begin{align}
    f_k(x) &=C\int_{\mathbb{R}^n}\frac{f_{k-1}(y)}{|x-y|^{n-\alpha}}dy+C_{k-1},\\
    v_k^R(x) &= \int_{\mathbb{R}^n}G^2_R(x,y)f_k(y)dy,
\end{align}
for $k=2,3,...,m-1$. Use similar process and induction, one proves that
\begin{align}
    u_{m-k}(x) &= C\int_{\mathbb{R}^n}\frac{f_{k}(y)}{|x-y|^{n-\alpha}}dy\mbox{ and }C_k=0\quad\mbox{for }k=1,3,...m-1,\\
    u(x) &= C\int_{\mathbb{R}^n}\frac{f_{m}(y)}{|x-y|^{n-\alpha}}dy+C_m\quad\mbox{for some }C_m\geq0.
\end{align}
Apply property \ref{prop riesz} to $\{f_k\}$ one obtains that
\begin{equation}
    u(x)=C\int_{\mathbb{R}^n}\frac{u^p(y)}{|x-y|^{n-2m-\alpha}}dy+C_m.
\end{equation}
Note that $u\geq C_m$. If $C_m>0$, then
\begin{equation}
    +\infty>u(x)>C\int_{\mathbb{R}^n}\frac{u^p(y)}{|x-y|^{n-2m-\alpha}}dy\geq\int_{\mathbb{R}^n}\frac{C_m}{|x-y|^{n-2m-\alpha}}dy=+\infty.
\end{equation}
Thus $C_m=0$, and $u$ indeed solves integral equation (\ref{integral}).
\end{proof}

\bibliographystyle{siam}
\bibliography{ref}
\end{document}